\documentclass[10pt]{amsart}
\newlength{\basicwidth}
\newlength{\shortbasicwidth}
\newlength{\basicheight}

\numberwithin{equation}{section}

\begin{document}

\title{Inequalities for Taylor series involving the divisor function}
\maketitle

\vspace{1cm}
\begin{center}
HORST ALZER$^a$  \quad\mbox{and} \quad MAN KAM KWONG$^b$
\end{center}

\vspace{0.75cm}
\begin{center}
$^a$ Morsbacher Stra{\ss}e 10, 51545 Waldbr\"ol, Germany\\
\emph{Email:}  \tt{h.alzer@gmx.de}
\end{center}

\vspace{0.3cm}
\begin{center}
$^b$ Department of Applied Mathematics, The Hong Kong Polytechnic University,\\
Hunghom, Hong Kong\\
{\emph{Email:}} \tt{mankwong@connect.polyu.hk}
\end{center}

\vspace{2.8cm}
{\bf{Abstract.}}  Let 
$$
T(q)=\sum_{k=1}^\infty d(k) q^k, \quad |q|<1,
$$
where
 $d(k)$ denotes the number of positive divisors of the natural number $k$. We present monotonicity properties of functions defined in terms of $T$.
 More specifically, we proved that
$$  H(q) := T(q)- \frac{\log(1-q)}{\log(q)}  $$
 is strictly increasing in $ (0,1) $ while
$$  F(q) := \frac{1-q}{q} \,H(q)  $$
 is strictly decreasing in $ (0,1) $.
 
These results are then applied to obtain various inequalities,
one of which states
 that the double-inequality
$$
\alpha \,\frac{q}{1-q}+\frac{\log(1-q)}{\log(q)}
< T(q)<
\beta \,\frac{q}{1-q}+\frac{\log(1-q)}{\log(q)}, \quad 0<q<1,
$$
holds with the best possible constant factors $\alpha=\gamma$ and $\beta=1$. Here, $\gamma$ denotes Euler's constant. This refines a result of Salem, who proved the inequalities with $\alpha=1/2$ and $\beta=1$.

\vspace{0.6cm}

\vspace{0.1cm}

\newpage

\section{Introduction}

In this paper, we study the Taylor series
$$
T(q)=\sum_{k=1}^\infty d(k) \,q^k, \quad |q|<1,
$$
where  $d(k)$ denotes the number of positive divisors of the natural number $k$. It is well-known that  the function $T$ has a close connection to Lambert series. We have
\begin{equation}
T(q)=
\sum_{k=1}^\infty\frac{q^k}{1-q^k}, \quad |q|<1.
\end{equation}
A proof of (1.1) and further information on Lambert series can be found in Knopp \cite[section 58 C]{K}. 
Another series representation for $T$  was given by Clausen \cite{C} in 1828,
$$
T(q)=\sum_{k=1}^\infty \frac{1+q^k}{1-q^k}\,q^{k^2}, \quad |q|<1.
$$
In 1899, Landau \cite{L} proved that $T$ can be used to determine  the value of a  series involving
 the classical Fibonacci numbers, defined by $F_0=F_1=1$, $F_{n}=F_{n-1}+F_{n-2}$ $(n\geq 2)$, 
$$
\sum_{k=1}^\infty \frac{1}{F_{2k}}=\sqrt{5}\,\bigl( T(c)-T(c^2)\bigr) = 1.53537..., \quad c=\Bigl(\frac{\sqrt{5}-1}{2}\Bigr)^2.
$$
Stimulated by his work on the analysis of data structure, Uchimura \cite{U} presented in 1981 the following result,
$$
T(q)=   (q;q)_{\infty} \sum_{k=1}^\infty \frac{k q^k}{(q;q)_k}, \quad |q|<1,
$$
where $(a;q)_k$ is the $q$-shifted factorial,
$$
(a;q)_k=\prod_{j=0}^{k-1} (1-a q^j), \quad (a;q)_{\infty}=\prod_{j=0}^\infty (1-a q^j).
$$
A related result was given by Merca \cite{M1}. In 2015, he proved
$$
T(q)=\frac{1}{(q;q)_{\infty}} \sum_{k=1}^\infty (-1)^{k-1}\,\frac{k q^{k+1\choose 2}}{(q;q)_k}, \quad |q|<1,
$$
and one year later,  he showed that there is a  relationship between partitions and $T$,
$$
T(q)=\frac{1}{(q;q)_{\infty}} \sum_{k=1}^\infty \bigl( s_{o}(k)-s_{e}(k)\bigr)\,q^k, \quad |q|<1.
$$
Here, $s_{o}(k)$ and $s_{e}(k)$ denote the number of parts in all partitions of $k$ into odd, respectively even, number of distinct parts; see Merca \cite{M2}.

The  $q$-digamma function is
 the logarithmic derivative of the $q$-gamma function, $\psi_q=\Gamma'_q/\Gamma_q$. The properties of $\Gamma_q$ and $\psi_q$ were investigated by numerous authors. 
For detailed  information on these functions we refer to  Askey \cite{A}, Salem \cite{S1, S2}, Salem \&  Alzahrani  \cite{SA} and the references cited therein. In view of the series representation
$$
\psi_q(x)=-\log(1-q)+(\log(q))\sum_{k=1}^\infty \frac{q^{kx }}{1-q^k}, \quad 0<q<1, \, x>0,
$$
we conclude from (1.1) that  $T$ can be expressed in terms of  $\psi_q(1)$,
\begin{equation}
T(q)=\frac{\psi_q(1)+\log(1-q)}{\log(q)}.
\end{equation}
The work on this paper has been inspired by an interesting double-inequality discovered by  Salem \cite{S3}. He proved 
\begin{equation}
0<1-\frac{1-q}{q\log(q)}\,\psi_q(1)<\frac{1}{2}, \quad 0<q<1.
\end{equation}
Using (1.2) and (1.3) we obtain elegant upper and lower bounds for $T(q)$. We have
\begin{equation}
\alpha \,\frac{q}{1-q}+\frac{\log(1-q)}{\log(q)}
< T(q)<
\beta \,\frac{q}{1-q}+\frac{\log(1-q)}{\log(q)}, \quad 0<q<1,
\end{equation}
with $\alpha=1/2$ and $\beta=1$. It is natural to ask whether these inequalities can be refined.
 More precisely, we look for the largest number $\alpha$ and the smallest number $\beta$ such that (1.4) is valid. Here, we  solve this problem. It turns out that $\beta=1$ is the best possible constant on the right-hand side of (1.4), but the factor $1/2$ on the left-hand side can be replaced by a larger number, namely by Euler's constant  $\gamma=0.57721...$. This reveals a connection between the divisor function and ``the third number of holy trinity ($\pi$, $e$, $\gamma$) of mathematical constants" \cite[p. 302]{B}.

In the next section, we collect several lemmas. Monotonicity properties of the functions
\begin{equation}
H(q)=T(q)-\frac{\log(1-q)}{\log(q)} \quad\mbox{and} \quad F(q)=\frac{1-q}{q}\Bigl( T(q)-\frac{\log(1-q)}{\log(q)}\Bigr)
\end{equation}
are given in Section 3. Finally, in Section 4, we apply the monotonicity of $F$ to prove (1.4) with $\alpha=\gamma$, $\beta=1$ and
we present sharp upper and lower bounds  for the  three Taylor series
\begin{equation}
\sum_{k=1}^\infty \bigl( d(k+1)-d(k)\bigr) \,q^k, \quad
\sum_{k=1}^\infty  \sum_{j=1}^k d(j) \,q^k, \quad
\sum_{k=2}^\infty  \sum_{j=1}^{k-1} \frac{d(j)}{k-j} \,q^k.
\end{equation}
The algebraical and numerical computations have been carried out  using the computer program MAPLE 13.

\vspace{0.3cm}
\section{Lemmas}

Throughout the paper we maintain the
 notations introduced in this section. The following nine lemmas play an important role in the proof of Theorem 3.2 given in Section 3.
We define for real numbers $q\in (0,1)$,  $x>0$ and integers $n\geq 1$,
$$
C_q(n)=\sum_{j=1}^n \sigma_q(j) \quad\mbox{and} \quad D_q(n)=\sum_{j=1}^n \rho_q(j)
$$
with
$$
 \sigma_q(j)=\int_j^{j+1} \phi_q(x)\,dx -\phi_q(j+1), \quad \rho_q(j)=\int_j^{j+1} \phi_q(x) \,dx -\frac{1}{2}\,\bigl( \phi_q(j)+\phi_q(j+1)\bigl)
$$
and
$$
\phi_q(x)=\frac{q^x(q^x -qx +x-1)}{(1-q^x)^2}.
$$
Geometrically, $ C_q(n) $ is the error of approximating the integral
$ \displaystyle \int_{1}^{n+1} \phi_q(x)\,dx $
using a special implementation of the rectangular rule.
Likewise $ D_q(n) $ is error of approximating the same integral using the
trapezoidal rule.

\vspace{0.3cm}
{\bf{Lemma 2.1.}}  \emph{Let $x \geq 1$. Then, $q\mapsto \phi_q(x)$ is increasing on $(0,1)$.}

\vspace{0.3cm}
\begin{proof}
We have for $q\in (0,1)$,
$$
\frac{\partial}{\partial q} \,\phi_q(x)=\frac{x^2(x^2-1)q^{x-1}}{(1-q^x)^3}\int_q^1 \int_t^1 (1-s) s^{x-2}\,ds\,dt\geq 0.
$$
\end{proof}

\vspace{0.1cm}
{\bf{Lemma 2.2.}} \emph{Let $q\in (0,1)$.}\\
\par\vspace*{-6mm}\par
\begin{itemize}
\addtolength{\itemsep}{-12pt}
\item[\rm(i)] \emph{There exists a number $N_q>1$ such that $\phi_q$ is strictly concave on $[1, N_q]$ and strictly convex on $[N_q,\infty)$.}\\
\item[\rm(ii)] \emph{There exists a number $M_q \in (1,N_q)$ such that $\phi_q$ is strictly increasing on $[1, M_q]$ and strictly decreasing on $[M_q,\infty)$.}
\end{itemize}

\vspace{0.3cm}
\begin{proof}
(i) Let $x\geq 1$. Differentiation gives
\begin{equation}
\phi''_q(x)=\frac{-q^x \log(q)}{(1-q^x)^4} \,a_q(x)
\end{equation}
with
$$
a_q(x)=(1-q)\Bigl( -x\bigl( 1+4q^x +q^{2x})\log(q)-2(1-q^{2x}) +\frac{(1-q^{2x})\log(q)}{1-q}\Bigl).
$$
We have
$$
a''_q(x)=4 q^x \log^2(q) (1-q) \,b_q(x)
$$
with
$$
b_q(x)=\frac{-\log(q)}{1-q}\Bigl( x(1-q)(1+q^x)+q^x\Bigr)+q^x-2.
$$
Using
$$
\frac{-\log(q)}{1-q}\geq \frac{2}{1+q}
\quad\mbox{and}
\quad
x(1-q)(1+q^x)-(1+q)(1-q^x)=x(x^2-1)\int_q^1 \int_t^1 (1-s)s^{x-2}\,ds\,dt\geq 0
$$
gives
$$
b_q(x)\geq \frac{2}{1+q} \Bigl((1+q)(1-q^x)+q^x\Bigr)+q^x-2=\frac{(1-q)\,q^x}{1+q}>0.
$$
It follows that $a_q$ is strictly convex on $[1,\infty)$. Since
$$
a_q(1)=-q\,(1-q)(3+q)\int_q^1 \frac{(1-t)(t^2+t+6)}{t^2(t+3)^2}\,dt<0 \quad\mbox{and}
\quad
\lim_{x\to\infty} a_q(x)=\infty,
$$
it follows  that there exists a number $N_q>1$ such that $a_q$ is negative on $(1,N_q)$ and positive on $(N_q,\infty)$. From (2.1) we conclude that  $\phi_q$ is strictly concave on $[1, N_q]$ and strictly convex on $[N_q,\infty)$.

(ii) We have
\begin{equation}
\phi'_q(x)=\frac{-q^x \log(q)}{(1-q^x)^2} \Bigl( \frac{-x(1-q)(1+q^x)}{1-q^x}+1-\frac{1-q}{\log(q)}\Bigr).
\end{equation}
It follows  that
$$
\phi'_q(1)=\frac{q}{(1-q)^2}\big(q \log(q)+1-q\big)>0
\quad\mbox{and} \quad \lim_{x\to\infty} \phi'_q(x)=0.
$$
Since $\phi'_q$ is strictly decreasing on $(1,N_q]$ and strictly increasing on $[N_q,\infty)$, there exists a number $M_q\in (1,N_q)$ such that $\phi'_q>0$ on $[1,M_q)$ and $\phi'_q<0$ on $(M_q,\infty)$. This implies that $\phi_q$ is strictly increasing on $[1, M_q]$ and strictly decreasing on $[M_q,\infty)$.
\end{proof}

\vspace{0.3cm}
{\bf{Lemma 2.3.}}  \emph{Let $q\in (0,1)$. If there exists an integer $m\geq 1$ such that $C_q(m)>0$, then $C_q(n)\geq C_q(m)$ for $n\geq m$.}

\vspace{0.3cm}
\begin{proof}
We claim that $M_q\leq m+1$. Suppose that this were false; then
$m+1<M_q$. An application of Lemma 2.2 (ii) gives for $j\in \{1,...,m\}$: $\sigma_q(j)<0$. Thus, $C_q(m)<0$, contradicting our assumption. 

Let $r\geq 1$ be an integer. Again we apply Lemma 2.2 (ii) and obtain
$$
C_q(m+r)-C_q(m+r-1)=\sigma_q(m+r)> 0.
$$
Hence, $C_q(m)<C_q(m+1) < C_q(m+2) <...$.
\end{proof}

\vspace{0.3cm}
{\bf{Lemma 2.4.}} (i) \emph{If $q\in (0, 0.117]$, then $C_q(1)>0$.}\\
(ii) \emph{If $q\in (0.117, 0.91]$, then $C_q(39)>0$.}

\vspace{0.3cm}
\begin{proof}
(i)  Let $q\in (0, 0.117]$. We have
\begin{equation}
C_q(1)=\int_1^2 \phi_q(x)\,dx -\phi_q(2)=\frac{U(q)}{(q+1)^2 \log^2(q)}
\end{equation}
with
$$
U(q)=-\big(1-q^2 +q\log(q)\big) \,q\log(q)-(q+1)^2 (q-1-\log(q))\log(1+q).
$$
Applying $\log(1+q)\leq q$ and $q-1-\log(q)>0$ gives
\begin{equation}
U(q) \geq -\big(1-q^2 +q\log(q)\big)\,q\log(q)-q(q+1)^2(q-1-\log(q))=q V(-\log(q)),
\end{equation}
where
$$
V(y)=1+\big(1-2y-y^2\big)\,\mbox{e}^{-y} -(1+2y) \,\mbox{e}^{-2y} -\mbox{e}^{-3y}.
$$
Since $0<q\leq 0.117$, we get $y=-\log(q)\geq -\log(0.117)=2.145...$. Using
$$
V'(y)=(y^2-3)\,\mbox{e}^{-y} +4\,y \,\mbox{e}^{-2y}+3 \,\mbox{e}^{-3y}>0
$$
yields
\begin{equation}
V(y)\geq V(-\log(0.117))=0.0022... .
\end{equation}
From (2.3), (2.4) and (2.5) we obtain $C_q(1)>0$.

(ii) Let  $q\in (0.117,0.91]$.
Applying Lemma 2.1 gives that
$$
C_q(39)=\int_1^{40} \phi_q(x)\,dx-\sum_{k=1}^{40} \phi_q(k)=W_1(q)-W_2(q)
$$
is the difference of two increasing functions. Let $0.117\leq r\leq q\leq s\leq 0.91$. Then,
$$
C_q(39)\geq W_1(r)-W_2(s)=W(r,s), \quad\mbox{say}.
$$
We set
$$
r_k=0.117+\frac{k}{10^3}, \quad r'_k=0.835+\frac{k}{2\cdot 10^4}, \quad r''_k=0.9+\frac{k}{5\cdot 10^5},
$$
$$
s_k=0.117+\frac{k+1}{10^3}, \quad s'_k=0.835+\frac{k+1}{2\cdot 10^4}, \quad s''_k=0.9+\frac{k+1}{5\cdot 10^5}.
$$
By direct computation we find
$$
W(r_k, s_k)>0 \,\, (k=0,1,...,717), \quad    W(r'_k,s'_k)>0 \,\,  (k=0,1,...,1299), \quad W(r''_k, s''_k)>0 \,\,  (k=0,1,...,4999).
$$
This yields $C_q(39)>0$.
\end{proof}

\vspace{0.3cm}
{\bf{Lemma 2.5.}} \emph{Let $q\in [0.91,1)$. Then, $N_q\geq 14$.}

\vspace{0.3cm}
\begin{proof}
In view of Lemma 2.2 (i) it suffices to show that $\phi''_q(14)< 0$. Using  (2.1) we conclude that we have to prove that
$a_q(14)< 0$, or, equivalently, $G(q)< 0$, where
$$
G(q)=-\log(q) \Delta(q)+2(1-q)(q^{28}-1)
$$
with
$$
\Delta(q)=13 -14q+56q^{14} -56 q^{15}+15 q^{28}-14 q^{29}.
$$
Next, we apply Sturm's theorem to determine the number of distinct roots of a polynomial in an interval; see van der Waerden \cite[section 79]{W}. We obtain that
$\Delta$ has precisely one zero on $[0.91,1]$. Since $\Delta(0.91)=1.76...$ and $\Delta(1)=0$, we conclude that $\Delta$ is positive on $[0,91,1)$.
Using this result and
$$
-\log(q)\leq 1-q +\frac{11}{20} \,(1-q)^2
$$
yields
$$
\frac{G(q)}{1-q}\leq \Bigl(1+\frac{11}{20}\,(1-q)\Bigr)\Delta(q)+ 2(q^{28}-1)=G_0(q), \quad\mbox{say}.
$$
An application of Sturm's theorem gives that $G_0$ has precisely one zero on $[0.91,1]$. We have $G_0(0.91)=-0.0028...$ and $G_0(1)=0$. It follows that $G_0$ and $G$ are negative on $[0.91,1)$.
\end{proof}

\vspace{0.3cm}
{\bf{Lemma 2.6.}} \emph{Let $q\in (0,1)$ and let $j\geq 1$ be an integer.}\\
\hspace*{0.0mm} (i) \emph{If $N_q\in [j,j+1]$ and $\phi'_q(x) \geq -\omega$ for $x\in [j, j+1]$ with $\omega\geq 0$, then $\rho_q(j)\geq - \omega/2$.}\\
(ii) \emph{If $\phi_q$ is convex on $[j,j+1]$, then}
\begin{equation}
\rho_q(j)\geq -\frac{1}{8} \bigl( \phi'_q(j+1)- \phi'_q(j) \bigr).
\end{equation}

\vspace{0.3cm}
\begin{proof}
(i) We consider two cases.\\
\underline{$\vphantom{y}$Case 1}. $M_q\leq j$. We have
\begin{equation}
\frac{1}{2}\,\omega \geq \frac{1}{2}\int_j^{j+1} (-\phi'_q(x))\,dx=\frac{\phi_q(j)-\phi_q(j+1)}{2}.
\end{equation}
Applying (2.7) and Lemma 2.2 (ii) gives
$$
\frac{1}{2}\,\omega+\rho_q(j)\geq \int_j^{j+1} \phi_q(x)\,dx-\phi_q(j+1)\geq 0.
$$
\underline{$\vphantom{y}$Case 2}. $j<M_q$.  We have $j < M_q < N_q \leq j+1$ and $\phi_q(j)\leq \phi_q(M_q)$, \,$\phi_q(j+1)\leq \phi_q(M_q)$. Let $\epsilon =M_q-j>0$ and $\delta=j+1-M_q>0$. Then, $\epsilon+\delta=1$ and
\begin{eqnarray}
\frac{\phi_q(j)+\phi_q(j+1)}{2} & = & 
 \epsilon  \,\frac{\phi_q(j)+\phi_q(j+1)}{2} +\delta \,\frac{\phi_q(j)+\phi_q(j+1)}{2}  \\ \nonumber
& \leq & \epsilon \,\frac{\phi_q(j)+\phi_q(M)}{2} +\delta \,\frac{\phi_q(M)+\phi_q(j+1)}{2}. \\ \nonumber
\end{eqnarray} 
Since $\phi_q$ is concave on $[j, M_q]$, we conclude from the Hermite-Hadamard inequality
\begin{equation}
\epsilon \,\frac{\phi_q(j)+\phi_q(M_q)}{2} \leq \int_j^{M_q} \phi_q(x) \,dx.
\end{equation}
 $\phi_q$ is decreasing on $[M_q, j+1]$ and $\phi'_q+\omega \geq 0$ on $[M_q, j+1]$. This implies
\begin{eqnarray}
-\frac{1}{2}\,\omega +\delta\,\frac{\phi_q(M_q)+\phi_q(j+1)}{2} & \leq &
-\frac{1}{2}\,\omega \,\delta^2 +\delta \,\frac{\phi_q(M_q)+\phi_q(j+1)}{2}  \\ \nonumber
& = & \frac{1}{2} \,\delta \Bigl( 2 \phi_q(j+1)-\int_{M_q}^{j+1} \bigl( \phi'_q(x)+\omega\bigr) \,dx\Bigr) \\ \nonumber
& \leq & \delta \,\phi_q(j+1) \\ \nonumber
& \leq & \int_{M_q}^{j+1} \phi_q(x) \,dx. \nonumber
\end{eqnarray}
Combining (2.8), (2.9) and (2.10) gives $\rho_q(j)\geq - \omega/2$.

(ii) We have
\begin{equation}
\frac{(b-a)^2}{8} \bigl( f'(b)-f'(a) \bigr) -(b-a) \frac{f(a)+f(b)}{2} +\int_a^b f(x)\,dx
\end{equation}
$$
=\int_a^{(a+b)/2} \Bigl(\frac{a+b}{2}-x\Bigl) \bigl( f'(x)-f'(a) \bigr) \,dx+
\int_{(a+b)/2}^b  \Bigl(x-\frac{a+b}{2}\Bigl) \bigl( f'(b)-f'(x) \bigr) \,dx.
$$
Applying  (2.11)
with $f=\phi_q$, $a=j$, $b=j+1$ gives (2.6).
\end{proof}

\vspace{0.3cm}
{\bf{Lemma 2.7.}} \emph{Let $q\in [0.91,1)$. The function
\begin{equation}
\Theta_q(x)=\frac{-q^x \log(q)}{(1-q^x)^2} \left( \frac {-x(1-q)(1+q^x)}{1-q^x}+q-\frac{1-q}{\log(q)}\right)
\end{equation}
 is increasing on $(0,\infty)$.}

\vspace{0.3cm}
\begin{proof}
Let $x>0$. We set $s=1-q^x$. Then, $s\in (0,1)$ and
$$
\Theta_q(x) = \eta_q(s)
$$
with
$$
\eta_q(s)=\frac{1-s}{s^3} \Bigl( (1-q)(2-s)\log(1-s) -s \bigl( q\log(q) +q-1\bigr) \Bigr).
$$
Using $q\log(q)/(1-q)>-1$ gives
\begin{eqnarray} \nonumber
\frac{s^4}{s(2-s)(1-q)} \eta'_q(s) &  = & \frac{(s^2-6s+6)\log(1-s)}{s(s-2)}+\frac{q\log(q)}{1-q}-2  \nonumber \\
& > &  \frac{(s^2-6s+6)\log(1-s)}{s(s-2)}-3  \nonumber \\
  & = &  \frac{s^2-6s+6}{s(2-s)} \int_0^s \frac{t^4}{(1-t)(t^2-6t+6)^2} \,dt  \nonumber \\
 & > & 0. \nonumber 
\end{eqnarray}
Since
$$
\Theta'_q(x)=-(1-s)\log(q) \,\eta'_q(s),
$$
we conclude that $\Theta'_q(x)>0$.
\end{proof}

\vspace{0.3cm}
{\bf{Lemma 2.8.}} \emph{Let $q\in [0.91,1)$ and $x\geq 1$. Then, $\phi'_q(x)\geq -0.035$.}

\vspace{0.3cm}
\begin{proof}
Applying Lemmas 2.2 (i), 2.5, 2.7 and (2.2), (2.12) leads to
\begin{equation}
\phi'_q(x)\,\geq \,\phi'_q(N_q)\,\geq  \,\Theta_q(N_q)\,\geq \,\Theta_q(14).
\end{equation}
Using $(1-q)/\log(q)\leq -q$ yields
\begin{equation}
-\Theta_q(14)\leq \frac{-q^{14}\log(q)}{(1-q^{14})^2}\left(\frac{14(1-q)(1+q^{14})}{1-q^{14}}-2q\right)=h_1(q)\bigl( h_2(q)+h_3(q)\bigr)
\end{equation}
with
$$
h_1(q)=\frac{-q^{14}\log(q)}{1-q^{14}}, \quad h_2(q)=\frac{1-q}{1-q^{14}},
\quad
h_3(q)=\frac{1}{1-q^{14}}\left( \frac{14(1-q)(1+q^{14})}{1-q^{14}}-q-1\right).
$$
From the integral representations
$$
h'_1(q)=\frac{q^{13}}{(1-q^{14})^2}\int_{q^{14}}^1 \frac{1-t}{t}\,dt,
\quad
h'_2(q)=\frac{-182}{(1-q^{14})^2}\int_q^1 (1-t) t^{12}\,dt,
$$
$$
h'_3(q)=\frac{-38220}{(1-q^{14})^3}\int_q^1 y^{12}\int_y^1 \int_t^1 (27-26s)s^{12} \,ds\,dt\,dy
$$
we conclude that $h_1$ is increasing and that $h_2$ and $h_3$ are decreasing on $(0,1)$. The functions $h_1$ and $h_2$ are positive on $(0,1)$ and since $\lim_{q\to 1} h_3(q)=0$, also $h_3$ is positive on $(0,1)$. Using 
$\lim_{q\to 1} h_1(q)=\frac1{14}$, we obtain for $q\in [0.91,1)$,
\begin{equation}
h_1(q)\bigl(h_2(q)+h_3(q)\bigr)\leq \frac{1}{14}\,\bigl( h_2(0.91)+h_3(0.91) \bigr)=0.034... .
\end{equation}
From (2.13), (2.14) and (2.15) we find $\phi'_q(x) \geq -0.035$ for $x\geq 1$.
\end{proof}

\vspace{0.3cm}
{\bf{Lemma 2.9.}} \emph{Let $q\in [0.91,1)$. Then, $D_q(10)>0.036$.}

\vspace{0.3cm}
\begin{proof}
We have
$$
D_q(10)-0.036=\left( \int_1^{11} \phi_q(x)\,dx-0.036\right) -\left(\sum_{k=1}^{10} \phi_q(k)+\frac{1}{2} \,\phi_q(11)\right)=J_1(q)-J_2(q).
$$
Applying Lemma 2.1 gives that $J_1$ and $J_2$ are increasing on $[0.91,1)$. Let $0.91\leq r\leq q\leq s\leq 1$. Then,
$$
D_q(10)-0.036\geq J_1(r)-J_2(s)=J(r,s), \quad\mbox{say}.
$$
We set
$$
r_k=0.91+\frac{k}{10^4}, \quad s_k=0.91+\frac{k+1}{10^4}
$$
and use
$$
J_2(1)=\lim_{q\to 1-} J_2(q)=\frac{208609}{55440}.
$$
Since $J(r_k,s_k)>0$ for $k=0,1,...,899$, we conclude that $D_q(10)-0.036>0$ for $q\in [0.91,1)$. 
\end{proof}

\vspace{0.3cm}
\section{Monotonicity theorems}

We prove monotonicity properties of the two functions defined in (1.5).

\vspace{0.3cm}
{\bf{Theorem 3.1.}} \emph{The function
$$
H(q)=T(q)-\frac{\log(1-q)}{\log(q)}
$$
is positive and strictly increasing on} $(0,1)$.

\vspace{0.3cm}
\begin{proof}
Let $0<q<1$. Using (1.1)  gives
\begin{equation}
q H'(q)=\sum_{k=1}^\infty k\frac{q^k}{(1-q^k)^2}+\frac{q\log(q)+(1-q)\log(1-q)}{(1-q) \log^2(q)}.
\end{equation}
Let
\begin{equation}
K_q(x)=\frac{x q^x}{(1-q^x)^2} \quad (x>0).
\end{equation}
Since
$$
K'_q(x)=-\frac{q^x (1+q^x)}{(1-q^x)^3 }\int_{q^x}^1\frac{y^2+1}{y(y+1)^2}\,dy<0,
$$
we conclude that $K_q$ is strictly decreasing on $(0,\infty)$, so that we get
\begin{equation}
\sum_{k=1}^\infty K_q(k)>\int_1^\infty K_q(x)\,dx=\frac{\log(1-q^x)}{\log^2(q)}+\left . \frac{x q^x}{(1-q^x)\log(q)} \right|^{x=\infty}_{x=1}
= -\frac{q\log(q)+(1-q)\log(1-q)}{(1-q)\log^2(q)}.
\end{equation}
From (3.1), (3.2) and (3.3) we obtain $H'(q)>0$. Thus, $H$ is strictly increasing on $(0,1)$ with $H(q)>\lim_{p\to 0}H(p)=0$ for $q\in (0,1)$.
\end{proof}

\vspace{0.3cm}
With the help of the results given in the previous section we are able to prove the following result.

\vspace{0.3cm}
{\bf{Theorem 3.2.}}  \emph{The function
$$
F(q)=\frac{1-q}{q} \, H(q)=\frac{1-q}{q}\left( T(q)-\frac{\log(1-q)}{\log(q)}\right)
$$
is strictly decreasing on $(0,1)$.}

\vspace{0.3cm}
\begin{proof}
Let $q\in (0,1)$. Then,
$$
q^2 F'(q)=\sum_{k=1}^\infty \phi_q(k)-A_q
$$
with
$$
A_q=\frac{-1}{\log^2(q)} \Bigl( \big(1-q+\log(q)\big) \log(1-q)+q \log(q) \Bigr).
$$
Let
$$
\Phi_q(x)=\frac{\bigl( q^{x+1} -\bigl(1+\log(q)\bigr)q^x +1-q+\log(q) \bigr) \log(1-q^x)+x q^x (1-q)\log(q)}{(1-q^x)\log^2(q)}.
$$
Since
$$
\Phi'_q(x)=\phi_q(x),  \quad \Phi_q(1) =-A_q\quad\mbox{and} \quad \lim_{x\to\infty} \Phi_q(x)=0,
$$
we obtain
$$
A_q=\int_1^{\infty} \phi_q(x) \,dx.
$$
It follows that $F'(q)<0$ is equivalent to
\begin{equation}
\sum_{k=1}^\infty \phi_q(k)<\int_1^\infty \phi_q(x) \,dx.
\end{equation}
To prove (3.4) we consider two cases.

\underline{$\vphantom{y}$Case 1}. $0<q\leq 0.91$. From Lemmas 2.3  and 2.4 we obtain
$$
C_q(n) \geq C_q(1)>0 \quad\mbox{for} \,\,  q\in(0, 0.117], \, n\geq 1
$$
and
$$
C_q(n)\geq C_q(39)>0  \quad\mbox{for} \,\,  q\in (0.117,0.91], \, n\geq 39.
$$
Thus, for $q\in (0, 0.91]$, 
$$
0<\lim_{n\to\infty} C_q(n)=\int_{1}^{\infty} \phi_q(x) \,dx -\sum_{k=1}^{\infty} \phi_q(k).
$$

\underline{$\vphantom{y}$Case 2}. $0.91<q<1$. Let $\tilde{N}_q$ be an integer such that $\tilde{N}_q < N_q\leq \tilde{N}_q+1$. From Lemma 2.5 we obtain $\tilde{N}_q\geq 13$. An application of Lemma 2.2 (i) and the Hermite-Hadamard inequality gives $\rho_q(j)\geq 0$ for $j=11,...,\tilde{N}_q-1$.
This result and Lemma 2.9 yield
\begin{equation}
D_q(\tilde{N}_q-1)=D_q(10)+\sum_{j=11}^{\tilde{N}_q-1} \rho_q(j) \geq D_q(10) > 0.036.
\end{equation}
From Lemmas 2.6 (i) and 2.8 we obtain
\begin{equation}
\rho_q(\tilde{N}_q)\geq -\frac{1}{2}\,\big( 0.035 \big).
\end{equation}
Next, we apply Lemma 2.6 (ii). Since $\phi_q$ is convex on $[N_q,\infty)$, we obtain for $j\geq \tilde{N}_q+1$,
$$
\rho_q(j)\geq -\frac{1}{8}\bigl( \phi'_q(j+1)-\phi'_q(j)\bigr).
$$
Using this inequality and Lemma 2.8 leads to
\begin{equation}
\sum_{j=\tilde{N}_q+1}^{\infty} \rho_q(j) \geq -\frac{1}{8} \sum_{j=\tilde{N}_q+1}^{\infty} \bigl( \phi'_q(j+1)-\phi'_q(j) \bigr)
=\frac{1}{8} \phi'_q( \tilde{N}_q+1) \geq -\frac{1}{8}\,\big( 0.035 \big).
\end{equation}
Combining (3.5), (3.6) and (3.7) gives
\begin{equation}
\sum_{j=1}^{\infty} \rho_q(j)=D_q(\tilde{N}_q-1)+\rho_q(\tilde{N}_q)+\sum_{j=\tilde{N}_q+1}^{\infty} \rho_q(j) >0.036-\frac{1}{2}\,\big( 0.035 \big)
-\frac{1}{8}\,\big( 0.035 \big) =0.014... .
\end{equation}
We have
\begin{equation}
\sum_{k=1}^m \rho_q(k)=\int_1^{m+1} \phi_q(x)\,dx-\sum_{k=1}^{m+1} \phi_q(k)+\frac{1}{2} \, \phi_q(m+1).
\end{equation}
Since $\lim_{x\to\infty} \phi_q(x)=0$, we conclude from (3.8) and (3.9) that (3.4) holds.
\end{proof}

\vspace{0.3cm}
An application of Theorems 3.1 and 3.2 leads to upper and lower bounds for the ratio $H(r)/H(s)$.

\vspace{0.3cm}
{\bf{Corollary 3.3.}} \emph{For all real numbers $r$ and $s$ with $0<r<s<1$ we have}
$$
\frac{r(1-s)}{s(1-r)}<\frac{H(r)}{H(s)}<1.
$$

\vspace{0.3cm}
\section{Inequalities}

We show that the monotonicity property of the function $F$ (defined in (1.5)) can be used  to obtain sharp upper and lower bounds for $T(q)$ and the Taylor series given in (1.6).  First, we present the best possible constant factors in double-inequality (1.4).

\vspace{0.3cm}
{\bf{Theorem 4.1.}} \emph{For all real numbers $q\in (0,1)$ we have
\begin{equation}
\alpha \,\frac{q}{1-q}+\frac{\log(1-q)}{\log(q)} \,<\, T(q) \,<\, 
\beta \,\frac{q}{1-q}+\frac{\log(1-q)}{\log(q)}
\end{equation}
with the best possible constant factors $\alpha=\gamma$ and $\beta=1$. }

\vspace{0.3cm}
\begin{proof}
The inequalities  (4.1) are equivalent to
\begin{equation}
\alpha<F(q)<\beta \quad (0<q<1).
\end{equation}
Since
$$
\lim_{q\to 0}\frac{(1-q)\log(1-q)}{q\log(q)}=0,
$$
we find
\begin{equation}
\lim_{q\to 0} F(q)=\lim_{q\to 0} \frac{1-q}{q} \,T(q)=d(1)=1.
\end{equation}
From (1.2) we obtain
$$
F(q)=\frac{1-q}{q\log(q)}\,\psi_q(1).
$$
We have
$$
\lim_{q\to 1} \frac{1-q}{q\log(q)}=-1
\quad\mbox{and} \quad
 \lim_{q\to 1} \psi_q(1)=\psi(1)=-\gamma;
$$
see Krattenthaler \& Srivastava \cite{KS}. It follows that
\begin{equation}
\lim_{q\to 1} F(q)=\gamma.
\end{equation}
Using the limit relations (4.3), (4.4) and Theorem 3.2 we conclude that (4.2) holds with the best possible bounds $\alpha=\gamma$ and $\beta =1$.
\end{proof}

\vspace{0.3cm}
Next, we offer inequalities for the Taylor series whose coefficients 
   are $d(k+1)-d(k)$ $(k=1,2,...)$. We mention an interesting property of this  difference which was discovered by Tur\'an \cite[p. 39]{MSC}. For each $c>0$ there exists a natural number $k$ such that $d(k+1)-d(k)>c$.

\vspace{0.3cm}
{\bf{Theorem 4.2.}} 
\emph{For all real numbers $q\in (0,1)$ we have
\begin{equation}
\alpha_0+\frac{(1-q)\log(1-q)}{q \log(q)}
< \sum_{k=1}^\infty \bigl( d(k+1)-d(k)\bigr) \,q^k<
\beta_0+\frac{(1-q) \log(1-q)}{q \log(q)}
\end{equation}
with the best possible constants $\alpha_0=\gamma-1$ and $\beta_0=0$.}

\vspace{0.3cm}
\begin{proof}
Let $q\in (0,1)$. We define
$$
F_0(q)= \sum_{k=1}^\infty \bigl( d(k+1)-d(k)\bigr) \,q^k-\frac{(1-q)\log(1-q)}{q\log(q)}.
$$
Since
$$
1+ \sum_{k=1}^\infty \bigl( d(k+1)-d(k)\bigr) \,q^k =\sum_{k=0}^\infty d(k+1)\,q^k -\sum_{k=1}^\infty d(k) \,q^k =\left(\frac{1}{q}-1\right) \,T(q),
$$
we get
$$
F_0(q)=F(q)-1.
$$
Applying  Theorem 3.2 and the limit relations
$$
\lim_{q\to 0} F_0(q)=0 \quad\mbox{and} \quad \lim_{q\to 1} F_0(q)=\gamma-1
$$
we obtain (4.5)  with the best possible  constants $\alpha_0=\gamma-1$ and $\beta_0=0$.
\end{proof}

\vspace{0.3cm}
The coefficients of the series given in the following theorem are
the partial sums of the divisor function which are related  to the floor function.
We have
$$
\sum_{k=1}^n d(k)=\sum_{k=1}^n [n/k],
$$
where
 $[x]$ denotes the  greatest integer less than or equal to $x$. 
These sums have a nice geometric interpretation. They give the exact number of lattice points in the area $x>0$, $y>0$, $xy\leq n$; see P\'olya \& Szeg\"o \cite[p. 131]{PS}.
The study of the average order of $d(k)$  dates back to Dirichlet and was continued by Hardy, Landau and others; see Apostol 
\cite[section 3.5]{Ap}.

\vspace{0.3cm}
{\bf{Theorem 4.3.}} 
\emph{For all real numbers $q\in (0,1)$ we have
\begin{equation}
{\lambda}\,\frac{q}{(1-q)^2}+\frac{\log(1-q)}{(1-q) \log(q)}
< \sum_{k=1}^\infty  \sum_{j=1}^k d(j) q^k<
{\mu} \,\frac{q}{(1-q)^2}+ \frac{\log(1-q)}{(1-q) \log(q)}
\end{equation}
with the best possible constant factors $\lambda=\gamma$ and $\mu=1$.}

\vspace{0.3cm}
\begin{proof}
Let $q\in (0,1)$ and
$$
L(q)=\frac{(1-q)^2}{q} \left( \sum_{k=1}^\infty  \sum_{j=1}^k d(j) \,q^k-\frac{\log(1-q)}{(1-q)\log(q)}\right).
$$
Since
$$
\frac{1}{1-q} \,T(q)=\sum_{k=1}^\infty  \sum_{j=1}^k d(j) \,q^k,
$$
we find  $L(q)=F(q)$. Applying Theorem 3.2, (4.3) and (4.4) leads to  (4.6)  with the best possible constant factors $\lambda=\gamma$ and $\mu=1$.
\end{proof}

\vspace{0.3cm}
We conclude the paper with a companion of (4.6).

\vspace{0.3cm}
{\bf{Theorem 4.4.}} 
\emph{For all real numbers $q\in (0,1)$ we have
\begin{equation}
\lambda_0\frac{q \log(1-q)}{1-q}-\frac{\log^2(1-q)}{\log(q)}
< \sum_{k=2}^\infty  \sum_{j=1}^{k-1} \frac{d(j)}{k-j} \, q^k<
\mu_0\frac{q\log(1-q)}{1-q}- \frac{\log^2(1-q)}{\log(q)}
\end{equation}
with the best possible constant factors $\lambda_0=-\gamma$ and $\mu_0=-1$.}

\vspace{0.3cm}
\begin{proof}
Let $q\in (0,1)$. Using
$$
-\log(1-q) T(q)=\sum_{k=2}^\infty  \sum_{j=1}^{k-1} \frac{d(j)}{k-j} \, q^k
$$
gives for
$$
L_0(q)=\frac{1-q}{q\log(1-q)} \left(  \sum_{k=2}^\infty  \sum_{j=1}^{k-1} \frac{d(j)}{k-j} \,q^k+\frac{\log^2(1-q)}{\log(q)}\right)
$$
the representation
$
L_0(q)=-F(q)$. From Theorem 3.2 and (4.3), (4.4) we conclude that (4.7) is valid with the smallest constant $\lambda_0=-\gamma$ and the largest constant $\mu_0=-1$.
\end{proof}

\end{document}